\documentstyle[11pt]{article}
\begin{document}
\newcommand{\p}{\parallel }
\makeatletter
\makeatother
\newtheorem{th}{Theorem}[section]
\newtheorem{lem}{Lemma}[section]
\newtheorem{de}{Definition}[section]
\newtheorem{rem}{Remark}[section]
\newtheorem{cor}{Corollary}[section]
\renewcommand{\theequation}{\thesection.\arabic {equation}}

\title{{\bf Differential Forms and the Wodzicki Residue for Manifolds with Boundary}
\thanks{partially supported by
 MOEC and the 973 project.}}
\author{   Yong Wang\\
{\scriptsize \it Nankai Institute of Mathematics}\\
{\scriptsize \it Tianjin 300071, P.R.China}\\
{\scriptsize \it E-mail: wangy581@nenu.edu.cn}}
\date{}
\maketitle

\noindent {\bf Abstract}~~ In [3], Connes found a conformal
invariant using Wodzicki's 1-density and computed it in the case
of 4-dimensional manifold without boundary. In [14], Ugalde
generalized the Connes' result to $n$-dimensional manifold without
boundary. In this paper, we generalize the results of [3] and [14]
to the case of manifolds
with boundary.\\
\noindent{\bf Subj. Class.:}\quad Noncommutative global analysis;
Noncommutative differential geometry.
\section{Introduction}
\quad In 1984, Wodzicki discovered a trace on the algebra
$\Psi_{\rm cl}(X)$ of all classical pseudodifferential operators
on a closed compact manifold $X$ in [15], which vanishes if the
order of the operator is less than $-{\rm {dim}}X$. It turns out
to be the unique trace on this algebra up to rescaling.

Wodzicki' residue has been applied to many branches of
mathematics. Especially, it plays a prominent role in
noncommutative geometry. In [4], Connes proved that Wodzicki's
residue coincided with Dixmier's trace on pseudodifferential
operators of order $-{\rm {dim}}X$. Wodzicki's residue also had
been used to derive an action for gravity in the framework of
noncommutative geometry in [6],[9],[10].

\quad In [3], for an even dimensional compact oriented conformal
manifold $X$ without boundary, Connes constructs a canonical
Fredholm module $(H,F)$. Here $H$ is the Hilbert space of square
integrable forms of middle dimension: $H=L^2(X,\wedge ^l_cT^*X)$
with $l=\frac{1}{2}{\rm dim}X$ and $F=2P-1$ where $P$ is the
orthogonal projection on the image of $d$. By Hodge decomposition
theorem, we observe that $F$ preserves the finite dimensional
space of harmonic forms $H^l$, and $F$ restricted to $H\ominus
H^l$ is given by $$F=\frac{d\delta-\delta d}{d\delta+\delta d}.$$
Using the equality :
$${\rm  Wres}(f_{0}[F,f_{1}][F,f_{2}])=\int
_{X}f_{0}\Omega_{n}(f_{1},f_{2}),\eqno(1.1)$$ where $f_{0}$,
$f_{1}$, $f_{2}\in C^{\infty}(X)$, Connes defined an $n$-form
$\Omega_{n}(f_{1},f_{2})$ which is uniquely determined, symmetric
in $f_1$ and $f_2$, and conformally invariant. In particular, in
the 4-dimensional case, this differential form was explicitly
computed in [3] by the conformal
deformation way.\\
\indent In [13], Ugalde presented the computations in the six
  dimensional case for a whole family of differential forms related
to $\Omega_{n}(f_1,f_2)$. In [14], he gave an explicit expression
of $\Omega_{n}(f_1,f_2)$
in the flat case and indicated the way of computation in the general case.\\
 \indent The purpose of this paper is to generalize these results to the case of manifolds with boundary. \\
\indent To do so, we find first that Wodzicki's residue in $(1.1)$
should be replaced by Wodzicki' residue for manifolds with
boundary. For a detailed introduction to the residue for manifolds
with boundary see [5], where Fedosov etc. defined a residue on
Boutet de Monvel's algebra and proved that it is a unique
continuous trace. For a good summary also see [11]. In addition,
Grubb and Schrohe got this residue through asymptotic expansions
in [8]. Subsequently we will use operator $\pi ^{+}F$ in Boutet de
Monvel's algebra instead of $F$ in $(1.1)$ ($\pi ^{+}F$ will be
introduced in
Section 2.1).\\
\indent Secondly, we will use the form pair
$(\Omega_{n,\pi^+S,X}(f_{1},f_{2}),\Omega_{n-1,\pi^+S,Y}(f_{1},f_{2}))$
instead of $\Omega _{n}(f_{1},f_{2})$ where $Y=\partial X$ and
$(1.1)$ turns into:\\ $~~\\ \widetilde{{\rm
Wres}}(\pi^+f_0[\pi^{+}F,\pi^+f_1][\pi^+F,\pi^+f_2])$\\
$$~~~~~ =\int _{X}f_0\Omega_{n,\pi^{+}F,X}(f_1,f_2)+\int _{Y}f_0|_Y\Omega
_{n-1,\pi^{+}F,Y}(f_1,f_2),\eqno(1.2)$$
 where $f_{0}, f_{1}, f_{2} \in C^{\infty}(X)$; $f_0|_Y$ denotes that the restriction of $f_0$ on Y. Here
 $f_{0}$ is
 assumed to be
 independent of $x_{n}$ near the boundary, where $x'=(x_1,...,x_{n-1})$ are coordinates on $\partial X$
and $x_n$ is  the normal coordinate (In what follows,
$x_n$ always denotes the normal coordinate.).\\
\indent In Section 2, we briefly recall Boutet de Monvel's
calculus and Wodzicki's residue for manifolds with boundary.

\indent In Section 3, for a pseudodifferential operator $S$ of
order 0 with the transmission property acting on sections of a
vector bundle $E$ over $\widetilde X=X\cup_Y X$, we consider
$\widetilde {{\rm
Wres}}({\pi^{+}f_0}[{\pi}^{+}S,{\pi}^{+}{f_{1}}][{\pi}^{+}S,{\pi}^{+}{f_{2}}])$
where $\pi ^{+}S: C^\infty (X,E|_X)\rightarrow C^\infty (X,E|_X)$
is defined in Section 2.1. We also show that\\
$$
\Omega_{n,\pi^+S,X}(f_1,f_2)=\Omega_{n,S,\tilde{X}}(\bar{f_1},\bar{f_2})|_X~;\eqno(1.3)$$
 $~~\widetilde{{\rm
 Wres}}(\pi^+f_0[\pi^{+}S,\pi^+f_1][\pi^+S,\pi^+f_2])$\\
$$~~~~~=\int _{X}f_0\Omega_{n,\pi^{+}S,X}(f_1,f_2)+\int _{Y}f_0|_Y\Omega
_{n-1,\pi^{+}S,Y}(f_1,f_2).\eqno(1.4)$$ \noindent determine a
unique form pair  $(\Omega
_{n,\pi^+S,X}(f_1,f_2),\Omega_{n-1,\pi^+S,Y}(f_1,f_2))$ which is
symmetric in $f_1$ and $f_2$, where $\bar{f_{1}}, \bar{f_{2}}$ are
extensions on $\tilde X$ of $f_{1}, f_{2}$ and $f_{0}$ is
independent of $x_{n}$ near the boundary and $\Omega _{n,S,\tilde
X}(\bar{f_{1}},\bar {f_{2}})=\Omega _{n}(\bar{f_{1}},\bar
{f_{2}})$ is defined in [3]. Moreover,$\widetilde{{\rm
Wres}}({\pi^+f_0}[\pi^+S,\pi^+f_1][\pi^+S,\pi^+f_2])$
 is a
Hochschild 2-cocycle (see Section 3) over $C^\infty (X)$.\\
\indent In Section 4, for a Riemannian manifold $(X,g)$ which has
the product metric near the boundary, $(\tilde X,\tilde g)$ is the
associated double Riemannian manifold. When dim$X$ is even, then
$\Omega_{n-1,\pi^+F,Y}(f_{1}, f_{2})=0$ and we get the formula:
$$\widetilde{{\rm Wres}}({\pi^+f_0}[\pi^+F,\pi^+f_1][\pi^+F,\pi^+f_2])=\int_{X}f_0\Omega
_{n,\pi^+F,X}(f_1,f_2).\eqno(1.5)                 $$ So we define
subconformal manifolds and $\Omega_{n,\pi^+F,X}(f_{1},f_{2})$ is a
obvious subconformal invariant. When dim$X$ is odd and
$(S,E)=(F,H)$, where $H=L^{2}(X,\wedge
^{\frac{n+1}{2}}_{c}T^{\star}X)$ and $F$ is defined as before,
then $\Omega_{n,\pi^+F,X}(f_{1},f_{2})=0$. So we get:
$$\widetilde{{\rm Wres}}(\pi^+f_0[\pi^+F,\pi^+f_1][\pi^+F,\pi^+f_2])
=\int_{Y}f_0|_Y\Omega_{n-1,\pi^+F,Y}(f_1,f_2).\eqno
 (1.6)$$

\indent Subsequently, in Sections 5,6, we compute the expression
of $\Omega_{n-1,\pi^+F,Y}(f_1,f_{2})$ and get its explicit
expression for flat manifolds in the $x_{n}$-independent and the
$x_{n}$-dependent cases. In Section 7, when $n=3$, using the
normal coordinate
 way we prove the formula:
$$ \Omega_{2,\pi^+F,Y}(f_1,f_2)=\frac{3}{8}\pi\Omega_{2}(f_1|Y,f_2|Y)-6\pi^2\partial_{x_n}f_1|_{x_n=0}
\partial_{x_n}f_2|_{x_n=0}{\rm Vol}_Y.\eqno(1.7)$$

\noindent So
$$\Omega_{2,\pi^+F,Y}(f_1,f_2)+6\pi^2\partial_{x_n}f_1|_{x_n=0}
\partial_{x_n}f_2|_{x_n=0}{\rm
Vol}_Y\eqno(1.8)$$ may be considered as a conformal invariant of
$(X,g)$. The above results generalize [3] and [14]
to the case of manifolds with boundary.\\
\indent ~~For the rest of this paper, We will briefly use
$(\Omega_n(f_1,f_2),\Omega_{n-1}(f_1,f_2))$ instead of
$(\Omega_{n,{\pi}^{+}S,X}(f_1,f_2),\Omega_{n-1,{\pi}^{+}S,Y}(f_1,f_2))$
in this section. We refer $\Omega_n(f_1,f_2)$ in this paper (in
[3]) if $f_1$ and $f_2$ are functions on manifolds with (without)
boundary.
\section{Boutet de Monvel's Calculus and Residue for
Manifolds with Boundary  }

\quad In this section, we recall some basic facts about Boutet de
Monvel's calculus which
 will be used in the following. For more details, see [1], [7], [11]
 and [12]. \\

\noindent{\bf 2.1~~Boutet de Monvel's Algebra}\\

\indent Let $$ F:L^2({\bf R}_t)\rightarrow L^2({\bf
R}_v);~F(u)(v)=\int e^{-ivt}u(t)dt$$ denote the Fourier
transformation and $\Phi(\overline{{\bf R}^+}) =r^+\Phi({\bf R})$
(similarly define $\Phi(\overline{{\bf R}^-}$)), where $\Phi({\bf
R})$ denotes the Schwartz space and $$r^{+}:C^\infty ({\bf
R})\rightarrow C^\infty (\overline{{\bf R}^+});~ f\rightarrow
f|\overline{{\bf R}^+};~ \overline{{\bf R}^+}=\{x\geq0;x\in {\bf
R}\}.$$ We define $H^+=F(\Phi(\overline{{\bf R}^+}));~
H^-_0=F(\Phi(\overline{{\bf R}^-}))$ which are orthogonal to each
other. We have the following property:
 $h\in H^+~(H^-_0)$ iff $h\in C^\infty({\bf R})$
which has an analytic extension to the lower (upper) complex
half-plane $\{{\rm Im}\xi<0\}~(\{{\rm Im}\xi>0\})$ such that for
all nonnegative integer $l$, $$\frac{d^{l}h}{d\xi^l}(\xi)\sim
\sum^{\infty}_{k=1}\frac{d^l}{d\xi^l}(\frac{c_k}{\xi^k})$$
 as $|\xi|\rightarrow +\infty,{\rm Im}\xi\leq0~({\rm Im}\xi\geq0)$.\\
\indent Let $H'$ be the space of all polynomials and
$H^-=H^-_0\bigoplus H';~H=H^+\bigoplus H^-.$ Denote by
$\pi^+~(\pi^-)$
 respectively the projection on $H^+~(H^-)$. For
calculations, we take $H=\widetilde H=\{$rational functions having
no poles on the real axis$\}$ ($\tilde{H}$ is a dense set in the
topology of $H$). Then on
$\tilde{H}$, \\
$$\pi^+h(\xi_0)=\frac{1}{2\pi i}\lim_{u\rightarrow
0^{-}}\int_{\Gamma^+}\frac{h(\xi)}{\xi_0+iu-\xi}d\xi,\eqno (2.1)$$
where $\Gamma^+$ is a Jordan close curve included ${\rm Im}\xi>0$
surrounding all the singularities of $h$ in the upper half-plane
and $\xi_0\in {\bf R}$.
 Similarly, define $\pi^{'}$ on $\tilde{H}$,\\
$$\pi'h=\frac{1}{2\pi}\int_{\Gamma^+}h(\xi)d\xi.\eqno(2.2)$$
So, $\pi'(H^-)=0$. For $h\in H\bigcap L^1(R)$,
$\pi'h=\frac{1}{2\pi}\int_{R}h(v)dv$ and for $h\in
H^+\bigcap L^1(R)$, $\pi'h=0$.\\
\indent An operator of order $m\in {\bf Z}$ and type $d$ is a matrix\\
$$A=\left(\begin{array}{lcr}
  \pi^+P+G  & K  \\
   T  &  S
\end{array}\right):
\begin{array}{cc}
\   C^{\infty}(X,E_1)\\
 \   \bigoplus\\
 \   C^{\infty}(\partial{X},F_1)
\end{array}
\longrightarrow
\begin{array}{cc}
\   C^{\infty}(X,E_2)\\
\   \bigoplus\\
 \   C^{\infty}(\partial{X},F_2)
\end{array}.
$$
where $X$ is a manifold with boundary $\partial X$ and
$E_1,E_2~(F_1,F_2)$ are vector bundles over $X~(\partial X
)$.~Here,~$P:C^{\infty}_0(\Omega,\overline {E_1})\rightarrow
C^{\infty}(\Omega,\overline {E_2})$ is a classical
pseudodifferential operator of order $m$ on $\Omega$, where
$\Omega$ is an open neighborhood of $X$ and
$\overline{E_i}|X=E_i~(i=1,2)$. $P$ has an extension:
$~{\cal{E'}}(\Omega,\overline {E_1})\rightarrow
{\cal{D'}}(\Omega,\overline {E_2})$, where
${\cal{E'}}(\Omega,\overline {E_1})~({\cal{D'}}(\Omega,\overline
{E_2}))$ is the dual space of $C^{\infty}(\Omega,\overline
{E_1})~(C^{\infty}_0(\Omega,\overline {E_2}))$. Let
$e^+:C^{\infty}(X,{E_1})\rightarrow{\cal{E'}}(\Omega,\overline
{E_1})$ denote extension by zero from $X$ to $\Omega$ and
$r^+:{\cal{D'}}(\Omega,\overline{E_2})\rightarrow
{\cal{D'}}(\Omega, {E_2})$ denote the restriction from $\Omega$ to
$X$, then define
$$\pi^+P=r^+Pe^+:C^{\infty}(X,{E_1})\rightarrow {\cal{D'}}(\Omega,
{E_2}).\eqno(2.3)$$ In addition, $P$ is supposed to have the
transmission property; this means that, for all $j,k,\alpha$, the
homogeneous component $p_j$ of order $j$ in the asymptotic
expansion of the
symbol $p$ of $P$ in local coordinates near the boundary satisfies:\\
$$\partial^k_{x_n}\partial^\alpha_{\xi'}p_j(x',0,0,+1)=
(-1)^{j-|\alpha|}\partial^k_{x_n}\partial^\alpha_{\xi'}p_j(x',0,0,-1),$$
then $\pi^+P:C^{\infty}(X,{E_1})\rightarrow C^{\infty}(X,{E_2})$
by [12]. Let $G$,$T$ be respectively the singular Green operator
and the trace operator of order $m$ and type $d$. $K$ is a
potential operator and $S$ is a classical pseudodifferential
operator of order $m$ along the boundary (For detailed definition,
see [11]). Denote by $B^{m,d}$ the collection of all operators of
order $m$
and type $d$,  and $\mathcal{B}$ is the union over all $m$ and $d$.\\
\indent Recall $B^{m,d}$ is a Fr\'{e}chet space. The composition
of the above operator matrices yields a continuous map:
$B^{m,d}\times B^{m',d'}\rightarrow B^{m+m',{\rm max}\{
m'+d,d'\}}.$ Write $$A=\left(\begin{array}{lcr}
 \pi^+P+G  & K \\
 T  &  S
\end{array}\right)
\in B^{m,d},
 A'=\left(\begin{array}{lcr}
\pi^+P'+G'  & K'  \\
 T'  &  S'
\end{array} \right)
\in B^{m',d'}.$$ The composition $AA'$ is obtained by
multiplication of the matrices(For more details see [12]). For
example $\pi^+P\circ G'$ and $G\circ G'$ are singular Green
operators of type $d'$ and
$$\pi^+P\circ\pi^+P'=\pi^+(PP')+L(P,P').\eqno(2.4)$$ Here $PP'$ is the usual
composition of pseudodifferential operators and $L(P,P')$ called
leftover term is a singular Green operator of type $m'+d$. The
composition formulas of the above operator symbols will be
given in the following.\\

\noindent{\bf 2.2~~Noncommutative Residue for Manifolds with Boundary}\\

\indent We assume that $E_1=E_2=E$;$~F_1=F_2=F$ and
$b(x',\xi',\xi_n,\eta_n)$ is the symbol of a singular Green
operator $G$ (about the definitions of symbols, see [11, p.11]),
then $${\rm tr}(b)=\frac{1}{2\pi}\int
_{\Gamma^+}b(x',\xi',\xi_n,\xi_n)d\xi_n=\bar{b}(x',\xi')\eqno(2.5)$$
is a symbol on $Y$ and $\bar{b}_{1-n}$ is obtained from $b_{-n}$
(see [5]). Let ${\bf S}~({\bf S}')$ be the unit sphere about
$\xi~(\xi')$ and $\sigma(\xi)~(\sigma(\xi'))$ be the corresponding
canonical
$n-1~(n-2)$ volume form. Now we recall the main theorem in [5],\\
\noindent{\bf Theorem (Fedosov-Golse-Leichtnam-Schrohe)}~~{\it Let
$X$ and $\partial X$ be connected, ${\rm dim}X=n\geq3$,
$A=\left(\begin{array}{lcr
}\pi^+P+G &   K \\
T &  S    \end{array}\right)$ $\in B$, and denote by $p$, $b$ and
$s$
the local symbols of $P,G$ and $S$ respectively. Define:\\
$~~{\rm{\widetilde{Wres}}}(A)=\int_X\int_{\bf
S}{\rm{tr}}_E\left[p_{-n}(x,\xi)\right]\sigma(\xi)dx$
$$~~~~~+2\pi\int_ {\partial
X}\int_{\bf S'}\left\{{\rm
tr}_E\left[({\rm{tr}}b_{-n})(x',\xi')\right]+{\rm{tr}}
_F\left[s_{1-n}(x',\xi')\right]\right\}\sigma(\xi')dx',\eqno (2.6)
$$
\noindent Then~~ a) ${\rm \widetilde{Wres}}([A,B])=0 $, for any
$A,B\in\mathcal{B}$;~~ b) It is a unique continuous trace on
$\mathcal{B}/\mathcal{B}^{-\infty}$.\\}

\section{Properties of
$(\Omega_{n,\pi^{+}S,X}(f_1,f_2),\Omega_{n-1,\pi^{+}S,Y}(f_1,f_2))$}

\quad  Let $X$ be a compact $n$-dimensional manifold with boundary
$Y$ and $\widetilde X=X\bigcup_YX$. For a pseudodifferential
operator $S$ of order 0 with the transmission property acting on
the sections of a vector bundle $E$ over $\tilde{X}$, we
consider the composition:\\
$$\widetilde{P}=\left(\begin{array}{cc}
\pi^+{f_0} & 0  \\
 0  & 0
\end{array}\right)
\left[\left(\begin{array}{cc}
\pi^+S & 0  \\
0  & 0
\end{array}\right),\left(\begin{array}{cc}
\pi^+f_1 & 0  \\
0 & 0
\end{array}\right)\right]\left[\left(\begin{array}{cc}
\pi^+S & 0 \\
0 & 0
\end{array}\right),\left(\begin{array}{cc}
\pi^+
f_2& 0  \\
0  & 0
\end{array}\right)\right]$$
$~~~~~~ :=\pi^+f_0[\pi^+S,\pi^+f_1][\pi^+S,\pi^+f_2].
           $\\

\noindent with $f_0, f_1,f_2\in C^\infty(X)$ which is the set
$\{f|_X|f\in C^\infty(\tilde X)\}$. By Section 2, $\pi^+S:
C^\infty (X,E|_X)\rightarrow C^\infty (X,E|_X)$ is well defined
and $\pi^{+}f_i: C^\infty (X,E|_X)\rightarrow C^ \infty (X,E|_X)$
is just the multiplication by
$f_i$ for $i=0,1,2$ and $\widetilde{P}=\pi^+(f_0[S,f_1][S,f_2])+G$
where $G$ is some singular Green operator. By (2.6),\\
$$\widetilde{{\rm Wres}}(\widetilde P)=\int_{X}f_0{\rm wres}
[S,\bar{f_{1}}][S,\bar{f_{2}}]|_X+2\pi\int_Y{\rm wres}_{x'} {\rm
tr}(b).\eqno(3.1)$$ Here $\bar{f_1},\bar{f_2}$ are the extensions
on $\widetilde{X}$ of $f_1,f_2$ and $${\rm
wres}[S,\bar{f_{1}}][S,\bar{f_{2}}]=\int_{\bf S}{\rm
tr}_Ep_{-n}(x,\xi)\sigma(\xi)dx;~ {\rm
 wres}_{x'}{\rm tr}(b)=\int_{\bf S'}{\rm tr}_E({\rm
 tr}b_{-n})(x',\xi')\sigma(\xi')dx',
 \eqno(3.2)$$
 \noindent where $p_{-n},b_{-n}$ are respectively the order $-n$ symbols
 of $[S,\bar{f_{1}}][S,\bar{f_{2}}]$ and $G$.
 Write:
$$\Omega_{n}(f_1,f_2)={\rm wres}[S,\bar{f_1}][S,\bar{f_2}]|_X=\Omega_n(\bar{f_1},\bar{f_2})
|_X;\eqno(3.3)$$
$$f_0|_Y\Omega_{n-1}(f_1,f_2)=2\pi  {\rm
wres}_{x'}{\rm tr}(b)\eqno(3.4),$$ then we have
$$\widetilde{{\rm Wres}}(\pi^+f_0[\pi^{+}S,\pi^+f_1][\pi^+S,\pi^+f_2])
=\int _{X}f_0\Omega_{n}(f_1,f_2)+\int
_{Y}f_0|_Y\Omega_{n-1}(f_1,f_2).\eqno(3.5)$$
   By [14], we have:
$\Omega_{n}(f_1,f_2)=$ $$\int_{|\xi|=1}{\rm tr}\left[
\sum\frac{1}{\alpha'!\alpha''!\beta!\delta!}D^{\beta}_x(\bar{f_1})
D^{\alpha''+\delta}_x(\bar{f_2})\right.\times
\left.\partial^{\alpha'+\alpha''+\beta}_{\xi}(\sigma^S_{-j})
\partial^{\delta}_{\xi}D^{\alpha'}_x(\sigma^S_{-k})\right]\sigma(\xi)d^nx\left|_X\right.
,\eqno(3.6)
$$
\noindent where $\sigma_{-j}^S$ denotes the order $-j$ symbol of
$S$; $D^\beta_x=(-i)^{|\beta|}\partial^\beta_x$ and the sum is
taken over $|\alpha'|+|\alpha''|+|\beta|+|\delta|+j+k=n;
|\beta|\geq1,|\delta|\geq1; \alpha',\alpha'',\beta,\delta \in {\bf
Z}^n_+; j,k\in {\bf Z}_+.$ By (3.6), this is a global $n$-form
which is
independent of  the extensions of $f_1,f_2$.\\
\indent Subsequently, we discuss the existence and uniqueness of
$\Omega_{n-1}(f_1,f_2)$.\\
\indent Recall, for example see [5, p.26], if $A_1$ and $A_2$ are
pseudodifferential operators with the transmission property, then
the Green operator
$$G=L(A_2,A_1)=\pi^+A_2\circ\pi^+A_1-\pi^+(A_2\circ A_1)$$ has a
symbol $b_{a_2,a_1}.$ If $A_1$ and $A_2$ have symbols
$a_1(\eta_n)=a_1(x',x_n,\xi',\eta_n)$ and
$a_2(\xi_n)=a_2(x',x_n,\xi',\xi_n)$ respectively, then
$b_{a_2,a_1}$ has an asymptotic expansion
formula:\\
$$b_{a_2,a_1}(x',\xi',\xi_n,\eta_n)\sim
\sum_{j,l,m=0}^{\infty}\frac{(-1)^mi^{j+l+m}}{j!l!m!}
\partial^j_{\xi_n}\partial^m_{\eta_n}b_{\partial^j_{x_n}a_2|x_n=0,
\partial^l_{\eta_n}\partial^{l+m}_{x_n}a_1|x_n=0}(x',\xi',\xi_n,\eta_n).\eqno(3.7)$$
\noindent When $a_1,a_2$ are independent of $x_n$ near the
boundary
, then we have:\\
$$b_{a_2,a_1}(x',\xi',\xi_n,\eta_n)=\frac{1}{2\pi}\int_{\Gamma^+}\frac{a^+_2(v)-a^+_2(\xi_n)}{v-\xi_n}\circ'
\frac{a^-_1(\eta_n)-a^-_1(v)}{\eta_n-v}dv,\eqno(3.8)$$ \noindent
where $a^+_i(v)=\pi^+_va_i(x',0,\xi',v),
a^-_i(v)=\pi^-_va_i(x',0,\xi',v),  i=1,2$ and\\
$$f(x',\xi',\xi_n)\circ'
g(x',\xi',\eta_n)=\sum_{|\alpha|\geq0}\frac{(-i)^{|\alpha|}}{\alpha!}\partial^{\alpha}_{\xi'}f\partial^{\alpha}_{x'}
g.\eqno(3.9)$$ \noindent Since $\pi_v^+f(x)=0$ and
$\pi_v^-f(x)=f(x)$, we get
 $b_{a_2,a_1}=0$ if $a_1$ or $a_2=f(x)$ by (3.7) and (3.8). So
 $b_{\sigma(S),f_1}=0$ and $[\pi^+S,\pi^+f]=\pi^+[S,f],$ then
\begin{eqnarray*}
\pi^+f_0[\pi^+S,\pi^+f_1][\pi^+S,\pi^+f_2]&=&\pi^+f_0\pi^+[S,f_1]\pi^+[S,f_2]\\
&=&\pi^+f_0[\pi^+([S,f_1][S,f_2])+\pi'B]\\
&=&\pi^+(f_0[S,f_1][S,f_2])+\pi^+f_0\circ \pi'B,
\end{eqnarray*}
\noindent where $\pi'B=L([S,f_1],[S,f_2])$ (here we use $f_i$ instead of $\bar{f_i}$) whose symbol is $b$.\\
\indent
  In the following we assume that
$f_0$ is independent of $x_n$ near the boundary, then we have
$$\sigma_{-n}(\pi^+f_0\circ
\pi'B)=f_0(x',0)b_{-n}(x',\xi',\xi_n,\eta_n).\eqno(3.10)$$ We can
see it in the boundary chart by the equality (see [11, p.11])
$$(\pi^+f_0\circ
\pi'B)u(x',x_n)=(2\pi)^{-n}\int
e^{ix\xi}\Pi'_{\eta_n}[f_0(x')b(x',\xi',\xi_n,\eta_n)
(e^+u)^{\wedge}(\xi',\eta_n)]d\xi.\eqno(3.11)$$  By definition:
\begin{eqnarray*}
2\pi  {\rm wres}_{x'}{\rm tr}(b)&=&\int_{|\xi'|=1}2\pi  {\rm tr}
 \left[{\rm tr}\sigma_{-n}(\pi^+f_0\circ
\pi'B)(x',\xi')\right]\sigma(\xi')d^{n-1}x'\\
&=&f_0(x',0)\int_{|\xi'|=1}\int_{\Gamma^+}{\rm tr}
b_{-n}(x',\xi',\xi_n,\xi_n)d\xi_n\sigma(\xi')d^{n-1}x'\\
&=&f_0|_Y\Omega_{n-1}(f_1,f_2),~~~~~~~~~~~~~~(3.12)
\end{eqnarray*}
\noindent then\\
$$\Omega_{n-1}(f_1,f_2)
=\int_{|\xi'|=1}\int_{\Gamma^+}{\rm tr}
b_{-n}(x',\xi',\xi_n,\xi_n)d\xi_n\sigma(\xi')d^{n-1}x'\eqno(3.13)$$
\noindent is an $(n-1)$-form over $Y$.\\
\noindent {\bf Theorem 3.1}~~{\it For the fixed $S$, the form pair
$(\Omega_n(f_1,f_2),\Omega_{n-1}(f_1,f_2))$ is uniquely determined
by (3.5) and (3.6).\\}
 \noindent{\it Proof.}~~ $\Omega_n(f_1,f_2)$ is uniquely
 determined by (3.6). We assume that $\Omega_{n-1}'(f_1,f_2)$ also satisfies (3.5), then
$$\int _{Y}f_0|_Y\Omega_{n-1}(f_1,f_2)=\int _{Y}f_0|_Y\Omega_{n-1}'(f_1,f_2)$$
for any $f_0|_Y\in C^{\infty}(Y)$. (In fact, using a cut function,
for any $g\in C^{\infty}(Y)$, there exists a function $f\in
C^{\infty}(X)$ such that $f|_Y=g$ and $f$ is independent of $x_n$
near the boundary.) So
$\Omega_{n-1}(f_1,f_2)=\Omega_{n-1}'(f_1,f_2)$.
\hfill$\Box$\\
\noindent {\bf Proposition 3.2}~~{\it $\widetilde{\rm
Wres}({\pi^+f_0}\pi^+[S,f_1]\pi^+[S,f_2])$
 is a
Hochschild 2-cocycle (for definition, see {\rm [6]}) over
$C^\infty (X)$.\\} \noindent{\it
Proof.}~~This proposition comes from the relations:\\
$[S,fh]=[S,f]h+f[S,h]$ ;
$~\pi^+(f_1[S,f_2])=\pi^+f_1\pi^+[S,f_2]$;
$\pi^+f_0\pi^+f_1=\pi^+(f_1f_0)$\\
 and the trace property of $\widetilde{{\rm Wres}}$.
\hfill$\Box$\\
\noindent{\bf Remark:}~~$\int _{X}f_0\Omega_{n}(f_1,f_2)$ and
$\int_{Y}f_0|_Y\Omega_{n-1}(f_1,f_2)$ are not Hochschild 2-cocycle
over $C^\infty (X)$.\\
\noindent {\bf Proposition 3.3}~~{\it $\Omega_n(f_1,f_2)$ and
$\Omega_{n-1}(f_1,f_2)$ are symmetric in $f_1$ and $f_2$.\\}
\noindent{\it Proof.}~~By [14],
$\Omega_n(\overline{f_1},\overline{f_2})$ is symmetric in
$\overline{f_1}$ and $\overline{f_2}$, so $\Omega_n(f_1,f_2)$ is
symmetric in $f_1$ and $f_2$. By the trace property of
$\widetilde{{\rm Wres}}$ and the commutativity of $C^{\infty}(X)$,
we note that:
$$\widetilde{{\rm Wres}}({\pi^+f_0}\pi^+[S,f_1]\pi^+[S,f_2]
=\widetilde{{\rm Wres}}({\pi^+f_0}\pi^+[S,f_2]\pi^+[S,f_1])$$
\noindent So $\Omega_{n-1}(f_1,f_2)$ is also symmetric in
$f_1,f_2$ by (3.5). \hfill$\Box$\\
\noindent{\bf Remark:}~~The condition $"fS^2=S^2f"$ in the theorem
2.7 of
[14] is not used here.\\
\indent In the following, we write the expression of
$\Omega_{n-1}(f_1,f_2)$ in detail. Let:\\
$$\overline{b_{a_1,a_2}}:={\rm tr}(b_{a_1,a_2})=\frac{1}{2\pi}\int
_{\Gamma^+}b_{a_1,a_2}(x',\xi',\xi_n,\xi_n)d\xi_n.\eqno(3.14)$$
\noindent By
[5, p.27], we have the formula:\\
$$\overline{b_{a_1,a_2}}=\sum_{j,k=0}^{\infty}\frac{(-i)^{j+k+1}}{(j+k+1)!}\pi'_{\xi_n}
\left[\partial^j_{x_n}\partial^k_{\xi_n}
a^+_1(x',0,\xi',\xi_n)\circ'
\partial^{j+1}_{\xi_n}\partial^k_{x_n}a^-_2(x',0,\xi',\xi_n)\right].\eqno(3.15)$$
\noindent Using (2.2),(3.13),(3.14) and (3.15), one obtains:\\
$$\Omega_{n-1}(f_1,f_2)=
\int_{|\xi'|=1}\int^{+\infty}_{-\infty}\left\{{\rm
trace}\sum^{\infty}_{j, k=0}\frac{(-i)^{j+k+1}}{(j+k+1)!}\right.$$
$$\left.\times\left[\partial^j_{x_n}\partial^k_{\xi_n}
a^+_1(x',0,\xi',\xi_n)\circ'\partial^{j+1}_{\xi_n}\partial^k_{x_n}a_2(x',0,\xi',\xi_n)\right]_{-n}\right\}
d\xi_n\sigma(\xi')d^{n-1}x'\eqno(3.16)$$ \noindent since the +
+ parts vanish after integration with respect to $\xi_n$ (see [5, p.23]).\\
\indent For $\pi'B=L([S,f_1],[S,f_2])$, by [14] lemma 2.2, then
for $i=1,2$, we have:
$$a_i=\sigma[S,f_i]=\sum_{k\geq1}\sigma_{-k}[S,f_i]=\sum_{k\geq1}\left[\sum^k_{|\beta|=1}\frac{1}{\beta!}D^{\beta}_x(f_i)
\partial^{\beta}_{\xi}(\sigma^S_{-(k-|\beta|)})\right].\eqno(3.17)$$
\noindent  By (3.9), then:

$$\left[\partial^j_{x_n}\partial^k_{\xi_n}
a^+_1(x',0,\xi',\xi_n)\circ'\partial^{j+1}_{\xi_n}\partial^k_{x_n}a_2(x',0,\xi',\xi_n)\right]_{-n}$$
\begin{eqnarray*}
&=&\left[\sum_{r,l}\partial^j_{x_n}\partial^k_{\xi_n}
a^+_{1(r)}(x',0,\xi',\xi_n)\circ'\partial^{j+1}_{\xi_n}\partial^k_{x_n}a_{2(l)}(x',0,\xi',\xi_n)\right]_{-n}\\
&=&\left[\sum_{r,l}\sum_{|\alpha|\geq0}\frac{(-i)^{|\alpha|}}{\alpha!}\partial^j_{x_n}
\partial^\alpha_{\xi'}\partial^k_{\xi_n}
a^+_{1(r)}(x',0,\xi',\xi_n)\times
\partial^\alpha_{x'}\partial^{j+1}_{\xi_n}\partial^k_{x_n}a_{2(l)}(x',0,\xi',\xi_n)\right]_{-n}\\
&=&\sum
\frac{(-i)^{|\alpha|}}{\alpha!}\partial^j_{x_n}\partial^\alpha_{\xi'}\partial^k_{\xi_n}
a^+_{1(r)}(x',0,\xi',\xi_n)\times
\partial^\alpha_{x'}\partial^{j+1}_{\xi_n}\partial^k_{x_n}a_{2(l)}(x',0,\xi',\xi_n)
\end{eqnarray*}
\noindent where the sum is taken over $
r-k-|\alpha|+l-j-1=-n,~~r,l\leq-1,~~|\alpha|\geq0$ for the fixed
$j,k$ and $a^+_{1(r)}~(a_{2(l)})$ denotes the order $r~(l)$ symbol
of $a^+_1~(a_2)$. Using:
$$
a^+_{1(r)}=\pi^+_{\xi_n}a_{1(r)}
=\pi^+_{\xi_n}\left[\sum_{|\beta|=1}^{-r}\frac{(-i)^{|\beta|}}{\beta!}\partial^\beta
_x(f_1)\partial^{\beta}_{\xi}(\sigma_{r+|\beta|}^S)\right]
=\sum_{|\beta|=1}^{-r}\frac{(-i)^{|\beta|}}{\beta!}\partial^\beta
_x(f_1)\pi^+_{\xi_n}\partial^{\beta}_{\xi}(\sigma^S_{r+|\beta|});$$
$$a_{2(l)}
=\sum_{|\delta|=1}^{-l}\frac{(-i)^{|\delta|}}{\delta!}\partial^\delta
_x(f_2)\partial^{\delta}_{\xi}(\sigma_{l+|\delta|}^S),~~~~~~~~~~~$$

\noindent we have:$$\left[\partial^j_{x_n}\partial^k_{\xi_n}
a^+_1(x',0,\xi',\xi_n)\circ'\partial^{j+1}_{\xi_n}\partial^k_{x_n}a_2(x',0,\xi',\xi_n)\right]_{-n}$$
$$=\sum
\sum^{-r}_{|\beta|=1}\sum^{-s}_{|\delta|=1}\frac{(-i)^{|\alpha|+|\beta|+|\delta|}}{\alpha!\beta!\delta!}
\partial^j_{x_n}\left[\partial^\beta_x(f_1)\partial^\alpha_{\xi'}\partial^k_{\xi_n}\pi^+_{\xi_n}\partial^\beta_{\xi}
(\sigma_{r+|\beta|}^S)\right]|_{x_n=0}\times$$
$$\partial^\alpha_{x'}\partial^k_{x_n}\left[\partial^\delta_x(f_2)
\partial^{j+1}_{\xi_n}\partial^\delta_\xi\sigma_{(l+|\delta|)}^S\right]|
_{x_n=0}\eqno(3.18)$$
\noindent with the sum $\sum$ as before. By (3.16) and (3.18), we get:\\
$$\Omega_{n-1}(f_1,f_2)=\sum^{\infty}_{j,k=0}\sum
\sum^{-r}_{|\beta|=1}\sum^{-l}_{|\delta|=1}
\frac{(-i)^{j+k+1+|\alpha|+|\beta|+|\delta|}}{\alpha!\beta!\delta!(j+k+1)!}$$
$$\times
\int_{|\xi'|=1} \int^{+\infty}_{-\infty}{\rm trace}\left\{
\partial^j_{x_n}\left[\partial^\beta_x(f_1)\partial^\alpha_{\xi'}\partial^k_{\xi_n}\pi^+_{\xi_n}\partial^\beta_{\xi}
(\sigma_{r+|\beta|}^S)\right]|_{x_n=0}\right.$$
$$\left.\times
\partial^\alpha_{x'}\partial^k_{x_n}\left[\partial^\delta_x(f_2)\partial^{j+1}_{\xi_n}
\partial^\delta_\xi
\sigma_{(l+|\delta|)}^S\right]|
_{x_n=0}\right\}d\xi_n\sigma(\xi')d^{n-1}x'\eqno(3.19)$$ \noindent
with the sum $\sum$ as (3.18).
\section{The Even Dimensional Case}

\quad Let $(X,g)$ be an even dimensional, compact, oriented,
Riemannian manifold with boundary $Y$ and product metric near the
boundary. $(\widetilde{X},\widetilde{g})$ is the associated double
manifold. Let $(E,S)=(H,F)$ associated to
$(\widetilde{X},\widetilde{g})$ introduced by Section 1. Let the dimension of $X$ be $n$.
 Since $\int_{|\xi|=1}$\{the product of odd
number of $\xi_i$ $\}\sigma(\xi)=0$, then we have\\
\noindent{\bf Lemma 4.1}~~{\it $\Omega_n(f_1,f_2)=0$ when $n$ is
odd and $\Omega_{n-1}(f_1,f_2)=0$ when $n$ is even.\\}
Since $n$ is even, by (3.5) and Lemma 4.1,we get:\\
$$\widetilde{{\rm Wres}}({\pi^+f_0}[\pi^+F,\pi^+f_1][\pi^+F,\pi^+f_2])=\int_{X}f_0\Omega
_{n}(f_1,f_2).\eqno(4.1)
$$
\noindent {\bf Definition 4.2}~~ A subconformal manifold is an
equivalence of Riemannian manifolds. Two metrics $g$ and
$\widetilde{g}$ are said to be equivalent if $\widetilde{g}=e^\eta
g$, where $\eta$ satisfies $\star)$ condition i.e. $\eta\in
C^{\infty}(X)$;~ $\eta\bigcup\eta\in C^{\infty}({\widetilde{X}})$
where $\eta\bigcup\eta=\eta$ on both copies of $X$.\\
\noindent{\bf Example:}~~1)~$X={\bf R}^n_+$ and $f(x)$ is an even
function about $x_n$, take $f|{\bf R}^n_+=\eta$, then
$\eta\bigcup\eta$ satisfies
$\star)$ condition.\\
~~ 2)~$f(x)$ is independent of $x_n$ near the boundary.\\
~~3)~$f(x)\in C^{\infty}(X),f(x)=e^{\frac{1}{x^2_n-1}}f(x')$ near
the boundary and if not, $f(x)=0$.\\
\indent Since the smoothness of $\eta\bigcup\eta$ just depends on
a neighborhood of the boundary, so we get:\\
\noindent{\bf Proposition 4.3}~~{\it $\eta\in C^{\infty}(X)$
satisfies $\star)$ condition iff $\exists f\in
C^{\infty}(\widetilde {X})$ such that
$f|_{Y\times(-1,1)}=\eta\bigcup\eta|_{Y\times(-1,1)}$.\\}
\noindent{\bf Proposition 4.4}~~{\it $\Omega_n(f_1,f_2)$ is
subconformally invariant for the above subconformal manifold.\\}
\noindent{\it Proof:}\quad Let $\widetilde{g}=e^\eta g$, where
$\eta$ satisfies $\star)$ condition, so
$\widetilde{g}\bigcup\widetilde{g}=e^{\eta\bigcup\eta}g\bigcup g$
and $\eta\bigcup\eta\in C^{\infty}(\widetilde{X}).$ By [3] or [14]
$\Omega_n(\bar{f_1},\bar{f_2})$ is conformal invariant, then
$\Omega_{n,g\bigcup
g}(\bar{f_1},\bar{f_2})=\Omega_{n,\widetilde{g}\bigcup
\widetilde{g}}(\bar{f_1},\bar{f_2})$ and
$\Omega_{n,g}(f_1,f_2)=\Omega_{n,\widetilde{g}}(f_1,f_2)$, where
$\Omega_{n,g}(f_1,f_2)$ denotes $\Omega_{n}(f_1,f_2)$ associated
to $g$.
\hfill$\Box$\\
\indent By [2, p.339], we have\\
\noindent{\bf Theorem 4.5}~~{\it Let $[(X,g)]$ be a
$4$-dimensional subconformal manifold with boundary as in the
definition 4.2 and $[(\widetilde{X},\widetilde{g})]$ be the
associated subconformal manifold without boundary, then
$$\Omega_{4}(f_1,f_2)=\frac{1}{16\pi^2}\left[\frac{1}{3}r\langle
d\widetilde{f_1},d\widetilde{f_2}\rangle-\triangle\langle
d\widetilde{f_1},d\widetilde{f_2}\rangle+\langle \nabla
d\widetilde{f_1},\nabla
d\widetilde{f_2}\rangle-\frac{1}{2}(\triangle\widetilde{f_1})
(\triangle\widetilde{f_2})\right] {\rm Vol}|_X,\eqno(4.2)$$ where
$\widetilde{f_1},\widetilde{f_2}\in C^{\infty}({\widetilde{X}})$
are the extensions of $f_1,f_2$, $r$ the scalar curvature, Vol the
volume form on $\widetilde{X}$, $\triangle$ the Laplacian and
$\nabla$ the Levi-civita connection associated to any metric of
$[(\widetilde{X},\widetilde{g})]$.\\}

\section{$\Omega_{n-1}(f_1,f_2)$ for Flat Manifolds in the $x_{n}$-Independent Case}

\quad In the rest of this paper, $(X,g)$ always denotes an odd
dimensional, compact, oriented Riemannian manifold with boundary
$Y$ and product metric near the boundary. Similar to Section 4, we
let
$(E,S)=(L^2(\wedge^{\frac{n+1}{2}}_c{T^\star\widetilde{X}}),F)$,
then $\Omega_n(f_1,f_2)=0$ by Lemma 4.1. So for $f_0$
independent of $x_n$ near the boundary, we have\\
$$\widetilde{{\rm Wres}}(\pi^+f_0[\pi^+F,\pi^+f_1][\pi^+F,\pi^+f_2])
=\int_{Y}(f_0|_Y)\Omega_{n-1}(f_1,f_2).\eqno
 (5.1)$$
\indent In this section, we assume that $X$ is flat and $f_1,f_2$
are independent of $x_n$ near the boundary and write
$\Omega_{n-1,{\rm flat}}(f_1,f_2)$ instead of $\Omega_{n-1}(f_1,f_2)$.\\
\indent We follow the method in Section 4 in [14]. Since $X$ is
flat, so is $(\widetilde{X},\widetilde{g})$. Then by Proposition
3.1 in [14], we have $\sigma(F)=\sigma_L(F)$ is independent of $x$
where $\sigma(F)~(\sigma_L(F))$ is the symbol (leading symbol) of
$F$. Using this information we deduce from (3.19) $j=k=0$ and
$|\beta|=-r, |\delta|=-l.$ Let $\beta=(\beta',\beta''),
\delta=(\delta',\delta'')$, where $\beta',\delta'\in
Z_+^{n-1},~\beta'',\delta''\in Z_+$, then by $f_1,f_2$ are
independent of $x_n$ near the boundary, we have
$\beta''=\delta''=0$
and\\
$$\Omega_{n-1,{\rm
flat}}(f_1,f_2)=\sum\sum_{|\beta'|=-r}\sum_{|\delta|=-l}
\frac{(-i)^{1+|\alpha|-r-s}}{\alpha!\beta'!\delta'!}\partial^{\beta'}_{x'}f_1(x',0)\partial^{\alpha+\delta'}_{x'}f_2(x',0)\times$$
$$\int_{|\xi'|=1} \int^{+\infty}_{-\infty}{\rm trace}\left[
\pi^+_{\xi_n}\partial^{\alpha+\beta'}_{\xi'}\sigma_L(F)\times
\partial_{\xi_n}\partial^{\delta'}_{\xi'}\sigma_L(F)\right]
d\xi_n\sigma(\xi')d^{n-1}x',$$ \noindent where the sum is taken
over $r+s-|\alpha|-1=-n, r\leq -1,~ s\leq -1,~|\alpha|\geq 0$. We
get \noindent{\bf Lemma 5.1}~~{\it $$\Omega_{n-1,{\rm
flat}}(f_1,f_2)=\sum
\frac{(-i)^n}{\alpha!\beta'!\delta'!}\partial^{\beta'}_{x'}f_1(x',0)\partial^{\alpha+\delta'}_{x'}f_2(x',0)\times$$
$$\int_{|\xi'|=1} \int^{+\infty}_{-\infty}{\rm trace}\left[
\pi^+_{\xi_n}\partial^{\alpha+\beta'}_{\xi'}\sigma_L(F)\times
\partial_{\xi_n}\partial^{\delta'}_{\xi'}\sigma_L(F)\right]
d\xi_n\sigma(\xi')d^{n-1}x',\eqno(5.2)$$ \noindent where the sum
is taken over $|\beta'|+|\delta'|+|\alpha|=n-1, |\beta'|\geq1,
|\delta'|\geq1$.\\}

\indent To better handle the previous expression, we consider:\\
$$\phi(\xi',\xi_n,u,v):=\sum\frac{1}
{\alpha!\beta'!\delta'!}u^{\beta'}v^{\alpha+\delta'}{\rm
trace}\left[
\pi^+_{\xi_n}\partial^{\alpha+\beta'}_{\xi'}\sigma_L(F)\times
\partial^{\delta'}_{\xi'}\partial_{\xi_n}\sigma_L(F)\right]\eqno(5.3)$$

\noindent with the sum as before and $u,v\in {\bf R}^{n-1}$. Then
by a recursive way we
have:\\
$$\Omega_{n-1,{\rm flat}}(f_1,f_2)=(-i)^n\left[\sum
A_{a,b}\partial^a_{x'}f_1(x',0)\partial^b_{x'}f_2(x',0)\right]d^{n-1}x'$$
\noindent where $A_{ab}$ is a number satisfying $\sum
A_{a,b}u^av^b=\int_{|\xi'|=1}
\int^{+\infty}_{-\infty}\phi(\xi',\xi_n,u,v)d\xi_n\sigma(\xi')$
and the sum is taken over  $a+b=n-1$ and $a\geq 1,~b\geq 1,~a,b\in
Z_{n-1}^+.$
 Instead of a direct approach to compute ${\rm trace}\left[
\pi^+_{\xi_n}\partial^{\alpha+\beta'}_{\xi'}\sigma_L(F)\times
\partial^{\delta'}_{\xi'}\partial_{\xi_n}\sigma_L(F)\right]$, we shall use
the Taylor expansion of function:\\
$$\psi(\xi',\eta',\xi_n):={\rm trace}\left[
\pi^+_{\xi_n}\sigma_L(F)(\xi',\xi_n)\times
\partial_{\xi_n}\sigma_L(F)(\eta',\xi_n)\right].$$
\noindent Considering the Taylor expansion of
$\psi(\xi'+u,\eta'+v,\xi_n)$ about $u,v$ at $u=v=0$, then:\\
$$\psi(\xi'+u,\eta'+v,\xi_n)=\sum_{
 |\beta|\geq0}\sum
 _{ |\delta|\geq0}
\frac{u^\beta v^\delta}{\beta!\delta!}{\rm
trace}\left[\partial^{\beta}_{\xi'}\pi^+_{\xi_n}\sigma_L(F)(\xi',\xi_n)\times
\partial^\delta_{\eta'}\partial_{\xi_n}\sigma_L(F)(\eta',\xi_n)\right]$$
\noindent with
$(\beta,\delta)=(\alpha_1,...,\alpha_{n-1},\alpha_n,...,\alpha_{2(n-1)}).$\\
\indent Write
$\psi(\xi',\eta',u,v,\xi_n):=\psi(\xi'+u,\eta'+v,\xi_n)$ and
$$T'_{n-1}\psi(\xi',\eta',u,v,\xi_n):=\sum
\frac{u^\beta v^\delta}{\beta!\delta!}{\rm
trace}\left[\partial^{\beta}_{\xi'}\pi^+_{\xi_n}\sigma_L(F)(\xi',\xi_n)\times
\partial^\delta_{\eta'}\partial_{\xi_n}\sigma_L(F)(\eta',\xi_n)\right],\eqno(5.4)$$
\noindent where the sum is taken over $|\beta|+|\delta|=n-1,~
 |\beta|\geq1,~|\delta|\geq1$ i.e. term of order $n-1$ in
the Taylor expansion of $\psi(\xi'+u,\eta'+v,\xi_n)$ minus the
terms with only powers of $u$ or only powers of $v$. Now, write:
$$P={\rm trace}\left[\partial^{\beta}_{\xi'}\pi^+_{\xi_n}\sigma_L(F)(\xi',\xi_n)\times
\partial^\delta_{\eta'}\partial_{\xi_n}\sigma_L(F)(\eta',\xi_n)\right],$$
\noindent then:
$$
T'_{n-1}\psi(\xi',\eta',u+v,v,\xi_n)=\sum \frac{(u+v)^\beta
v^\delta}{\beta!\delta!}=\sum \sum_{\beta'+\beta''=\beta}
 \frac{u^{\beta'}
v^{\beta''+\delta}}{\beta'!\beta''!\delta!}P$$
$$=\sum
\sum_{
  \beta'+\beta''=\beta;
  \beta'\neq0
}
 \frac{u^{\beta'}
v^{\beta''+\delta}}{\beta'!\beta''!\delta!}P\\
+\sum
 \frac{
v^{\beta+\delta}}{\beta!\delta!}P$$ where the sum $\sum$ is taken
over $
 |\beta|+|\delta|=n-1;
~|\beta|\geq1; ~|\delta|\geq1.$
$$=\sum
 \frac{u^{\beta'}
v^{\beta''+\delta}}{\beta'!\beta''!\delta!}P+
T'_{n-1}\psi(\xi',\eta',v,v,\xi_n).$$ where the sum $\sum$ is
taken over $|\beta'|+|\delta'|+|\delta|=n-1;~
|\beta'|\geq1;~|\delta|\geq1.$
 Therefore, by taking
$\eta=\xi$ we obtain:\\
$$T'_{n-1}\psi(\xi',\xi',u+v,v,\xi_n)-T'_{n-1}\psi(\xi',\xi',v,v,\xi_n)
=\phi(\xi',\xi_n,u,v).$$
\noindent In summary, we have:\\
\noindent{\bf Theorem 5.2}~~{\it $$\Omega_{n-1,{\rm
flat}}(f_1,f_2)=(-i)^n\left[\sum
A_{a,b}\partial^a_{x'}f_1(x',0)\partial^b_{x'}f_2(x',0)\right]d^{n-1}x',\eqno(5.5)$$
\noindent where $\sum A_{a,b}u^av^b=\int_{|\xi'|=1}
\int^{+\infty}_{-\infty}\left[T'_{n-1}\psi(\xi',\xi',u+v,v,\xi_n)
-T'_{n-1}\psi(\xi',\xi',v,v,\xi_n)\right]d\xi_n\sigma(\xi').$ and
$T'_{n-1}\psi(\xi',\eta',u,v,\xi_n)$ is defined by (5.4).\\}
\indent By Theorem 5.2, to obtain an explicit expression of
$\Omega_{n-1}$ in the flat case, it is necessary to study
$\psi(\xi',\eta',\xi_n)$
 for
$\xi'$ and $\eta'$ not zero in $T^{\star}_xY$. Recall
 the theorem 4.3 in [14] (we will find it is also correct when $n$ is odd
 through the check.) says that: when
 $\sigma_L(F)(\xi)\sigma_L(F)(\eta)$ acts on $m$-forms on $\widetilde{X}$, then
 $${\rm
trace}\left[\sigma_L(F)(\xi)\times \sigma_L(F)(\eta)\right]=
a_{n,m}\frac{\langle\xi,\eta\rangle^2}{|\xi|^2|\eta|^2}+b_{n,m}\eqno(5.6)$$
\noindent where
$b_{n,m}=C_n^m-a_{n,m}=C^{m-2}_{n-2}+C^m_{n-2}-2C^{m-1}_{n-2}$ and
$C^m_n$ denotes a combinator number. Suppose that $g=g^{Y}+d^2x_n$
near the boundary, then
$$\langle\xi,\eta\rangle_g=\langle\xi',\eta'\rangle_{g^Y}+\xi_n\eta_n\eqno(5.7)$$
where $\xi=\xi'+\xi_ndx_n;~\eta=\eta'+\eta_ndx_n.$ By (5.6) and
(5.7),
 then
\begin{eqnarray*}
{\rm trace}[\pi^+_{\xi_n}\sigma_L(F)(\xi',\xi_n)\times
\partial_{\xi_n}\sigma_L(F)(\eta',\xi_n)]
&=&\pi^+_{\xi_n}\partial_{\eta_n}{\rm
trace}[\sigma_L(F)(\xi)\times
\sigma_L(F)(\eta)]|_{\eta_n=\xi_n}\\
&=&\pi^+_{\xi_n}\partial_{\eta_n}\left[a_{n,m}\frac{\langle\xi,\eta\rangle^2}{|\xi|^2|\eta|^2}+b_{n,m}\right]
\left|_{\eta_n=\xi_n}\right.\\
&=&a_{n,m}\pi^+_{\xi_n}\left[\frac{2\langle\xi,\eta\rangle\xi_n|\eta|^2-2\eta_n\langle\xi,\eta\rangle^2}{|\xi|^2|\eta|^4}
\right]|_{\eta_n=\xi_n},
\end{eqnarray*}
\noindent by (2.1), Cauchy integral formula and the choice of
$\Gamma^+$:
\begin{eqnarray*}
&=&\frac{a_{n,m}}{|\eta|^4}\frac{1}{2\pi i}{\rm lim}_{u\rightarrow
0^-}\int_{\Gamma^+}\frac
{2(\langle\xi',\eta'\rangle+z\eta_n)z|\eta|^2-2\eta_n(\langle\xi',\eta'\rangle+z\eta_n)^2}
{(|\xi'|^2+z^2)(\xi_n+iu-z)}dz|_{\eta_n=\xi_n}\\
&=&\frac{a_{n,m}}{|\eta|^4} \frac
{2(\langle\xi',\eta'\rangle+i|\xi'|\eta_n)i|\xi'\|\eta|^2-
2\eta_n(\langle\xi',\eta'\rangle+i|\xi'|\eta_n)^2}
{2i|\xi'|(\xi_n-i|\xi'|)}|_{\eta_n=\xi_n}\\
&=&\frac{a_{n,m}(\langle\xi',\eta'\rangle+i|\xi'|\xi_n)}
{i|\xi'|(\xi_n-i|\xi'|)(|\eta'|^2+\xi_n^2)^2}
\left[|\eta'|^2i|\xi'|-\xi_n\langle\xi',\eta'\rangle\right].
\end{eqnarray*}
\noindent So we have:\\
\noindent{\bf Theorem 5.3}~~Suppose that $(X,g)$ has a product
metric near the boundary. When
 $\sigma_L(F)(\xi',\xi_n)\sigma_L(F)(\eta',\xi_n)$ acting on $m$-forms in the
 boundary chart, then\\
 {\it $~~~~~~{\rm
trace}\left[\pi^+_{\xi_n}\sigma_L(F)(\xi',\xi_n)\times
\partial_{\xi_n}\sigma_L(F)(\eta',\xi_n)\right]$\\
$$~~~~~~=\frac{a_{n,m}(\langle\xi',\eta'\rangle+i|\xi'|\xi_n)}
{i|\xi'|(\xi_n-i|\xi'|)(|\eta'|^2+\xi_n^2)^2}
\left[|\eta'|^2i|\xi'|-\xi_n\langle\xi',\eta'\rangle\right],\eqno(5.8)$$
\noindent where
$a_{n,m}=C^m_n-C^{m-2}_{n-2}-C^m_{n-2}+2C^{m-1}_{n-2}$.\\}

\section{$\Omega_{n-1}(f_1,f_2)$ for Flat Manifolds in the $x_{n}$-Dependent Case}

\quad In this section, we assume that $X$ is flat and $f_1,f_2$
are
dependent of $x_n$ near the boundary.\\
\indent Since $X$ is flat, so $\sigma(F)=\sigma_L(F)$ and
$|\beta|=-r,|\delta|=-s.$   By (3.19), we have:\\
\noindent{\bf Lemma 6.1}~~{\it
$$\Omega_{n-1,{\rm flat}}(f_1,f_2)=
\sum^{\infty}_{j,k=0}
\sum\frac{(-i)^n}{\alpha!\beta'!\beta''!\delta'!\delta''!(j+k+1)!}
\partial^{j+\beta''}_{x_n}\partial^{\beta'}_{x'}f_1|_{x_n=0}\times
\partial^{\alpha+\delta'}_{x'}\partial^{k+\delta''}_{x_n}f_2|_{x_n=0}\times$$
$$ \int_{|\xi'|=1}\int^{+\infty}_{-\infty}
{\rm
trace}\left\{\partial^{\alpha+\beta'}_{\xi'}\partial^{\delta'}_{\eta'}
\left[\partial^{k}_{\xi_n}
\pi^+_{\xi_n}\partial^{\beta''}_{\xi_n}\sigma_L(F)(\xi',\xi_n)\times
\partial^{j+1+\delta''}_{\xi_n}\sigma_L(F)(\eta',\xi_n)\right]|_{\xi'=\eta'}\right\}
d\xi_n \sigma(\xi')d^{n-1}x',\eqno(6.1)$$ \noindent where the sum
is taken over
$|\beta'|+\beta''+|\delta'|+\delta''+|\alpha|+j+k+1=n,
~~|\beta'|+\beta''\geq1,~~ |\delta'|+\delta''\geq1,
~~|\alpha|\geq0.$\\}
 \indent Similar to Section 5,
we want to give its explicit expression. Let:\\
$$\widetilde{\phi}(\xi',\xi_n,u,v):=\sum\frac{1}{\alpha!\beta'!\beta''!\delta'!\delta''!(j+k+1)!}u^{(\beta',j+\beta'')}
v^{(\alpha+\delta',k+\delta'')} $$
$$\times{\rm trace}\left[\partial^{\alpha+\beta'}_{\xi'}\partial^{k}_{\xi_n}
\pi^+_{\xi_n}\partial^{\beta''}_{\xi_n}\sigma_L(F)(\xi',\xi_n)\times
\partial^{\delta'}_{\xi'}\partial^{j+1+\delta''}_{\xi_n}\sigma_L(F)(\xi',\xi_n)\right]
$$
\noindent with the sum as before and $u,v \in {\bf R}^n$. One obtains:\\
$$\Omega_{n-1,{\rm flat}}(f_1,f_2)=(-i)^n\left[\sum
A_{a,b}\partial^a_{x}f_1(x',0)\partial^b_{x}f_2(x',0)\right]d^{n-1}x'$$
\noindent with $\sum A_{a,b}u^av^b=\int_{|\xi'|=1}
\int^{+\infty}_{-\infty}\widetilde{\phi}(\xi',\xi_n,u,v)d\xi_n\sigma(\xi')$.\\
\noindent Now,\\
$$\widetilde{\phi}(\xi',\xi_n,u,v)=\sum^{\infty}_{j,k=0}\sum_{\beta'',\delta''}
\frac{u^{j+\beta''}_nv^{k+\delta''}_n}{\beta''!\delta''!(j+k+1)!}\widetilde{\phi}_{j,k,\beta'',\delta''}
(\xi',\eta',\xi_n,u',v')|_{\eta'=\xi'}\eqno(6.2)$$\noindent with
\\
$$\widetilde{\phi}_{j,k,\beta'',\delta''}
(\xi',\eta',\xi_n,u',v')=\sum\frac{1}{\alpha!\beta'!\delta'!}
{u'}^{\beta'}{v'}^{\alpha+\delta'}$$
$$\times\partial^{\alpha+\beta'}_{\xi'}\partial^{\delta'}_{\eta'}{\rm trace}\left[
\partial^{k}_{\xi_n}
\pi^+_{\xi_n}\partial^{\beta''}_{\xi_n}\sigma_L(F)(\xi',\xi_n)\times
\partial^{j+1+\delta''}_{\xi_n}\sigma_L(F)(\eta',\xi_n)\right]\eqno(6.3)$$
\noindent where the sum is taken over
$|\beta'|+|\delta'|+|\alpha|=n-(j+k+1)-\beta''-\delta''=s,
~~|\beta'|+\beta''\geq1,~~ |\delta'|+\delta''\geq1$ for fixed $j,k,\beta'',\delta''.$\\
 \noindent
Write\\ $$\psi_{j,k,\beta'',\delta''} (\xi',\eta',\xi_n):={\rm
trace}\left[\partial^{k}_{\xi_n}
\pi^+_{\xi_n}\partial^{\beta''}_{\xi_n}\sigma_L(F)(\xi',\xi_n)\times
\partial^{j+1+\delta''}_{\xi_n}\sigma_L(F)(\eta',\xi_n)\right];$$

$$\psi_{j,k,\beta'',\delta''}
(\xi',\eta',\xi_n,u',v'):=\psi_{j,k,\beta'',\delta''}
(\xi'+u',\eta'+v',\xi_n);$$
$$T_s\psi_{j,k,\beta'',\delta''}
(\xi',\eta',u',v',\xi_n):=\sum_{|\beta|+|\delta|=s}\frac{u'^\beta
v'^\delta}{\beta!\delta!}$$
$$\times{\rm trace}\left[\partial^{\beta}_{\xi'}\partial^{k}_{\xi_n}
\pi^+_{\xi_n}\partial^{\beta''}_{\xi_n}\sigma_L(F)(\xi',\xi_n)\times
\partial^{\delta}_{\eta'}\partial^{j+1+\delta''}_{\xi_n}\sigma_L(F)(\eta',\xi_n)\right]$$
\noindent i.e. the term of order $s$ in the Taylor expression  of
$\psi_{j,k,\beta'',\delta''}(\xi',\eta',u,v,\xi_n)$. By (6.3),\\
$$\widetilde{\phi}_{j,k,\beta''\neq0,\delta''\neq0}
(\xi',\eta',\xi_n,u',v')=\sum_{|\beta'|+|\delta'|+|\alpha|=s}\frac{1}{\alpha!\beta'!\delta'!}u'^{\beta'}v'^{\alpha+\delta'}
\times
\partial^{\alpha+\beta'}_{\xi'}\partial^{\delta'}_{\eta'}\psi_{\beta''\neq0,\delta''\neq0}(\xi',\eta',\xi_n)$$
\begin{eqnarray*}
T_s\psi_{\beta''\neq0,\delta''\neq0}(\xi',\eta',u'+v',v')&=&\sum_{|\beta|+|\delta|=s}\sum_{\beta'+\beta''=\beta}
\frac{u'^{\beta'}v'^{\beta''+\delta}}{\beta'!\beta''!\delta!}
\partial^{\beta'+\beta"}_{\xi'}\partial^{\delta}_{\eta'}
\psi_{j,k,\beta''\neq0,\delta''\neq0}\\
&=&\widetilde{\phi}_{\beta''\neq0,\delta''\neq0}(\xi',\eta',u',v',\xi_n).~~~~~~~~(6.4)
\end{eqnarray*}
\noindent Let $T_{s,u}(T_{s,v})$ denote the term of order $s$ in
the Tayler expansion $\psi_{j,k,\beta'',\delta''}$ minus the terms
with only powers of $v~(u)$ and $T_s'$ denote the term of order
$s$ in the Tayler expansion $\psi_{j,k,\beta'',\delta''}$ minus
the terms with only powers of $u$ or $v$. In a similar way, we
get:\\
$$\begin{array}{lcr
} \
\widetilde{\phi}_{\beta''=0,\delta''=0}(\xi',\eta',u',v',\xi_n)=T'_s\psi_{\beta''=0,\delta''=0}(\xi',\eta',u'+v',v')-
T'_{s}\psi_{\beta''=0,\delta''=0}(\xi',\eta',v',v'); \\

\
\widetilde{\phi}_{\beta''=0,\delta''\neq0}(\xi',\eta',u',v',\xi_n)=T_{s,u}\psi_{\beta''=0,\delta''\neq0}(\xi',\eta',u'+v',v')-
T_{s,u}\psi_{\beta''=0,\delta''\neq0}(\xi',\eta',v',v'); \\

\
\widetilde{\phi}_{\beta''\neq0,\delta''=0}(\xi',\eta',u',v',\xi_n)=T_{s,v}\psi_{\beta''\neq0,\delta''=0}(\xi',\eta',u'+v',v')-
T_{s,v}\psi_{\beta''\neq0,\delta''=0}(\xi',\eta',v',v').
\end{array}\eqno(6.5)$$
 \noindent Summarizing, we have a similar result for
manifolds with boundary to the
theorem 4.2 in [14]:\\
\noindent{\bf Theorem 6.2}~~{\it $$\Omega_{n-1,{\rm
flat}}(f_1,f_2)=(-i)^n\left[\sum
A_{a,b}\partial^a_{x}f_1(x',0)\partial^b_{x}f_2(x',0)\right]d^{n-1}x'$$
\noindent with $\sum A_{a,b}u^av^b=\int_{|\xi'|=1}
\int^{+\infty}_{-\infty}\widetilde{\phi}(\xi',\xi_n,u,v)d\xi_n\sigma(\xi')$
and $\widetilde{\phi}(\xi',\xi_n,u,v)$ is determined by (6.2)
(6.4) and (6.5).\\} The computation of
$\psi_{j,k,\beta'',\delta''}(\xi',\eta',\xi_n)$
is similar to the theorem 5.3.\\

\section{The Computation of $\Omega _{2}(f_1,f_2)$ when $n=3$}

\quad In this section, we will give the global expression of
$\Omega_{2}(f_1,f_2)$ in three cases.\\
\noindent{\bf a)~~Flat and $f_1,f_2$ Are Independent of $x_n$ Near
the
Boundary.}\\
\indent By lemma 5.1 and $n=3$, we have $|\delta'|=|\beta'|=1$, $|\alpha|=0$ and\\
$$\Omega_{2,{\rm flat}}(f_1,f_2)=\sum^2_{i,j=1}(-i)^3
\partial_{x_i}f_1(x',0)\partial_{x_j}f_2(x',0)$$
$$\times\int_{|\xi'|=1}\int^
 {+\infty}_{-\infty}\partial_{\xi_i}\partial_{\eta_j}\left\{
 {\rm trace}\left[\pi^+_{\xi_3}\sigma_L(F)\times
 \partial_{\xi_3}\sigma_L(F)\right]\right\}|_{\xi'=\eta'}d\xi_3\sigma(\xi')dx_1\wedge
 dx_2.$$
 \noindent In this subsection we denote $\left|{{\begin{array}{cc}
\ _{\xi'=\eta'} \\  \ _{|\xi'|=1}
\end{array}}}\right.$ by $|_{\star}$. Using the theorem 5.3, then for $n=3~,m=2$ we have:
 \begin{eqnarray*}
D_{ij}:&=&\partial_{\xi_i}\partial_{\eta_j}\left\{{\rm
trace}\left[\pi^+_{\xi_n}\sigma_L(F)(\xi',\xi_n)\times
\partial_{\xi_n}\sigma_L(F)(\eta',\xi_n)\right]\right\}|_{\star}\\
&=&\partial_{\xi_i}\left(\frac{a_{n,m}}{\xi_n-i|\xi'|}A\right)|_{\star}
=\left[\frac{ia_{n,m}\xi_i}{(\xi_n-i|\xi'|)^2|\xi'|}A
+\frac{a_{n,m}}{\xi_n-i|\xi'|}\partial_{\xi_i}A\right]|_{\star} ,
\end{eqnarray*}
\noindent where\\
$$A=\partial_{\eta_j}\left\{\frac{\langle\xi',\eta'\rangle+i|\xi'|\xi_n}{|\eta'|^2+\xi^2_n}
\left[1-\frac{\xi_n(\langle\xi',\eta'\rangle+i|\xi'|\xi_n)}
{i|\xi'|(|\eta'|^2+\xi^2_n)}\right ]\right\}=A_1-A_2$$ \noindent and\\
$$A_1=\frac{\xi_j(|\eta'|^2+\xi^2_n)-2\eta_j(\langle\xi',\eta'\rangle+i|\xi'|\xi_n)}{(|\eta'|^2+\xi^2_n)^2}
\left[1-\frac{\xi_n(\langle\xi',\eta'\rangle+i|\xi'|\xi_n)}{i|\xi'|(|\eta'|^2+\xi^2_n)}\right
]$$
$$A_2=\frac{\langle\xi',\eta'\rangle+i|\xi'|\xi_n}{|\eta'|^2+\xi^2_n}\frac{\xi_n}{i|\xi'|}
\frac{\xi_j(|\eta'|^2+\xi^2_n)-2\eta_j(\langle\xi',\eta'\rangle+i|\xi'|\xi_n)}{(|\eta'|^2+\xi^2_n)^2}.$$
\noindent Through the computation,\\
$$~~~~~\partial_{\xi_i}A_1|_{\star}=
\frac{\delta_{ij}(1+\xi_n^2)-2\xi_i\xi_j(1+i\xi_n)}{(1+\xi^2_n)^2}
\left[1-\frac{\xi_n(1+i\xi_n)}{i(1+\xi^2_n)}\right ];$$
$$A_1|_{\star}=\frac{-i\xi_j}{(1-i\xi_n)^2(\xi_n+i)};~~~
B_1:=a_{n,m}\partial_{\xi_i}\left(\frac{A_1}{\xi_n-i|\xi'|}\right
)|_{\star}
=\frac{ia_{n,m}}{(\xi_n+i)^2(\xi_n-i)^2}\left(\delta_{ij}-\frac{i\xi_i\xi_j}{\xi_n+i}\right);
$$ $$A_2|_{\star} =\frac{i\xi_j\xi_n}{(1-i\xi_n)^3};~~~
\partial_{\xi_i}A_2|_{\star}
=\frac{\xi_n}{i(1+i\xi_n)(1-i\xi_n)^3}[\delta_{ij}(1-i\xi_n)-2\xi_i\xi_j];$$
$$B_2:=a_{n,m}\partial_{\xi_i}\left(\frac{A_2}{\xi_n-i|\xi'|}\right
)|_{\star}
=\frac{a_{n,m}\xi_n}{(\xi_n+i)^2(\xi_n-i)^2}\left(\delta_{ij}
-\frac{i\xi_i\xi_j}{\xi_n+i}\right);$$
 $$D_{ij}=B_1-B_2
=\frac{a_{n,m}}{(\xi_n-i)(\xi_n+i)^2}\left(-\delta_{ij}
+\frac{i\xi_i\xi_j}{\xi_n+i}\right).$$
 \noindent Using the fact that
$\int_{|\xi'|=1}\xi_i\xi_j\sigma(\xi')=\pi\delta_{ij}$ and
$\int_{|\xi'|=1}\sigma(\xi')=2\pi$ where $|\xi'|=1$ is the unit
circle, we have $i=j$ and
\begin{eqnarray*}
\Omega_{2,{\rm flat}}(f_1,f_2)
&=&\sum^2_{j=1}(-i)^3\partial_{x_j}f_1(x',0)\partial_{x_j}f_2(x',0)
\int_{|\xi'|=1}\int^
 {+\infty}_{-\infty}D_{jj}d\xi_3\sigma(\xi')dx_1\wedge
 dx_2\\
&=&i\sum^2_{j=1}\partial_{x_j}f_1(x',0)\partial_{x_j}f_2(x',0)
\int_{|\xi'|=1}\int
 _{\Gamma^+}
a_{3,2}\frac{-1+\frac{i\xi^2_j}{\xi_n+i}}{(\xi_n-i)(\xi_n+i)^2}
 d\xi_3\sigma(\xi')dx_1\wedge
 dx_2\\
&=&i\sum^2_{j=1}\partial_{x_j}f_1(x',0)\partial_{x_j}f_2(x',0)
\frac{a_{3,2}\pi
i}{2}\left[-\frac{1}{2}\int_{|\xi'|=1}\xi^2_j\sigma(\xi')+\int_{|\xi'|=1}\sigma(\xi')\right]\\
&=&-3\pi^2\sum^2_{j=1}\partial_{x_j}f_1(x',0)\partial_{x_j}f_2(x',0)dx_1\wedge
 dx_2\\
&=&-3\pi^2d(f_1|_{Y})\wedge\star d(f_2|_{Y}),
\end{eqnarray*}
\noindent because $X$ is flat and $a_{3,2}=4$.\\
{\bf b)~~Flat and $f_1,f_2$ Are Dependent of $x_n$ Near the
Boundary.}\\
\indent Since $n=3$ and $|\beta|\geq1,|\delta|\geq1,$ so we have
$|\delta|=|\beta|=1,|\alpha|=j=k=0$. By Lemma 6.1, then:\\
$$\Omega_{2,{\rm flat}}(f_1,f_2)
=\sum_{|\beta|=1}\sum_{|\delta|=1}(-i)^3\partial_{x_n}^{\beta''}\partial
^{\beta'}_{x'}f_1(x',0)\partial_{x'}^{\delta'}\partial
^{\delta''}_{x_n}f_2(x',0)\times$$
$$
\int_{|\xi'|=1}\int^
 {+\infty}_{-\infty}\partial_{\xi'}^{\beta'}\partial_{\eta'}^{\delta'}\left\{
 {\rm trace}\left[\pi^+_{\xi_n}\partial^{\beta''}_{\xi_n}\sigma_L(F)(\xi',\xi_n)\times
 \partial_{\xi_n}^{1+\delta''}\sigma_L(F)(\eta',\xi_n)\right]\right\}|_{\xi'=\eta'}d\xi_n\sigma(\xi')dx_1\wedge
 dx_2$$
$$=D_1+D_2+D_3+D_4,$$
\noindent where\\
$$D_1=\sum_{|\beta'|=1}\sum_{|\delta'|=1}(-i)^3\partial
^{\beta'}_{x'}f_1(x',0)\partial_{x'}^{\delta'}f_2(x',0)\times$$
$$\int_{|\xi'|=1}\int^
 {+\infty}_{-\infty}\partial_{\xi'}^{\beta'}\partial_{\eta'}^{\delta'}\left\{
 {\rm trace}\left[\pi^+_{\xi_n}\sigma_L(F)(\xi',\xi_n)\times
 \partial_{\xi_n}\sigma_L(F)(\eta',\xi_n)\right]\right\}|_{\xi'=\eta'}d\xi_n\sigma(\xi')dx_1\wedge
 dx_2;$$
$$D_2=\sum_{|\beta'|=1}(-i)^3\partial
^{\beta'}_{x'}f_1(x',0)\partial_{x_n}f_2(x',0)\times$$
$$\int_{|\xi'|=1}\int^
 {+\infty}_{-\infty}\partial_{\xi'}^{\beta'}\left\{
 {\rm trace}\left[\pi^+_{\xi_n}\sigma_L(F)(\xi',\xi_n)\times
 \partial_{\xi_n}^{2}\sigma_L(F)(\eta',\xi_n)\right]\right\}|_{\xi'=\eta'}d\xi_n\sigma(\xi')dx_1\wedge
 dx_2;$$
$$D_3=\sum_{|\delta'|=1}(-i)^3\partial_{x_n}
f_1(x',0)\partial_{x'}^{\delta'}f_2(x',0)\times$$
$$\int_{|\xi'|=1}\int^
 {+\infty}_{-\infty}\partial_{\eta'}^{\delta'}\left\{
 {\rm trace}\left[\pi^+_{\xi_n}\partial_{\xi_n}\sigma_L(F)(\xi',\xi_n)\times
 \partial_{\xi_n}\sigma_L(F)(\eta',\xi_n)\right]\right\}|_{\xi'=\eta'}d\xi_n\sigma(\xi')dx_1\wedge
 dx_2;$$
$$D_4=(-i)^3\partial_{x_n}
f_1(x',0)\partial_{x_n}f_2(x',0)\times$$
$$\int_{|\xi'|=1}\int^
 {+\infty}_{-\infty}{\rm trace}\left[\pi^+_{\xi_n}\partial_{\xi_n}\sigma_L(F)(\xi',\xi_n)\times
 \partial_{\xi_n}^{2}\sigma_L(F)(\xi',\xi_n)\right]d\xi_n\sigma(\xi')dx_1\wedge
 dx_2.$$
Observation: $D_1= -3\pi^2d(f_1|_{Y})\wedge\star d(f_2|_{Y})$
 by case a). In order to compute $D_2$, we must compute $ {\rm
trace}[\pi^+_{\xi_n}\sigma_L(F)(\xi',\xi_n)\times
 \partial_{\xi_n}^{2}\sigma_L(F)(\eta',\xi_n)]. $
Instead of the way of Theorem 5.3, we compute
$\pi^+_{\xi_n}\sigma_L(F)(\xi',\xi_n)$ firstly. Let
$p(\xi',\xi_n)=\varepsilon(\xi)l(\xi)-l(\xi)\varepsilon(\xi)$ be a
polynomial with matrices as coefficients where $\varepsilon(\xi)$
and $l(\xi)$ are the exterior and interior multiplications
respectively, then
$$\sigma_L(F)=\frac{p(\xi',\xi_n)}{|\xi'|^2+\xi^2_n}$$ by
Proposition 3.1 in [14].
By (2.1), we have:\\
$$\pi^+_{\xi_n}\left[\frac{p(\xi',\xi_n)}{|\xi'|^2+\xi^2_n}\right]=\frac{p(\xi',i|\xi'|)}{2i|\xi'|(\xi_n-i|\xi'|)};$$
$$\partial^2_{\xi_n}\sigma_{L}(F)(\eta',\xi_n)=\frac{\partial^2_{\xi_n}p(\eta',\xi_n)}{|\eta'|^2+\xi_n^2}-
\frac{4\xi_n\partial_{\xi_n}p(\eta',\xi_n)}{(|\eta'|^2+\xi_n^2)^2}-\frac{2|\eta'|^2-6\xi_n^2}{(|\eta'|^2+\xi_n^2)^3}
p(\eta',\xi_n).$$ \noindent Using $|$ instead of taking
$\xi_n=i|\xi'|,~\eta_n=\xi_n$,
by Theorem 4.3 of [14] ($n$=odd case), then\\
$$T:={\rm trace}[p(\xi)\times
p(\eta)]=a_{n,m}\langle\xi,\eta\rangle^2+b_{n,m}|\xi|^2|\eta|^2$$
\noindent so, $$
T|=a_{n,m}\left[\langle\xi',\eta'\rangle+i|\xi'|\xi_n\right]^2;~~
\partial_{\eta_n}T|=2ia_{n,m}|\xi'|\left[\langle\xi',\eta'\rangle+i|\xi'|\xi_n\right];~~
\partial^2_{\eta_n}T|=-2a_{n,m}|\xi'|^2;$$
$${\rm trace}\left[\pi^+_{\xi_n}\sigma_L(F)(\xi',\xi_n)\times
 \partial_{\xi_n}^{2}\sigma_L(F)(\eta',\xi_n)\right]
=\frac{1}{2i|\xi'|(\xi_n-i|\xi'|)}\times
\left[\frac{1}{|\eta'|+\xi^2_n}\partial^2_{\eta_n}T|\right.$$
$$\left.-\frac{4\xi_n}{(|\eta'|^2+\xi^2_n)^2}
\partial_{\eta_n}T|-\frac{2|\eta'|^2-6\xi_n^2}{(|\eta'|^2+\xi_n^2)^3}T|\right],
\eqno(7.1)$$
 \noindent then compute the partial derivative
$\partial_{\xi_i}$ of (7.1) and take $\xi'=\eta'$ and $|\xi'|=1$,
we have the result has form $\xi_if(\xi_n)$. Using
$\int_{|\xi'|=1}\xi_i\sigma(\xi')=0$, so $D_2=0$. Similarly, we
have $D_3=0$. In order to compute $D_4$, we'll compute $${\rm
trace}[\pi^+_{\xi_n}\partial_{\xi_n}\sigma_L(F)(\xi',\xi_n)\times
 \partial_{\xi_n}^{2}\sigma_L(F)(\xi',\xi_n)]$$. Since
\begin{eqnarray*}
\pi^+_{\xi_n}\partial_{\xi_n}\sigma_L(F)(\xi',\xi_n)&=&\pi^+_{\xi_n}\left[
\frac{\partial_{\xi_n}p(\xi',\xi_n)}{|\xi'|^2+\xi^2_n}-\frac{2\xi_np(\xi',\xi_n)}{(|\xi'|^2+\xi_n^2)^2}\right];\\
\pi^+_{\xi_n}\left[\frac{\partial_{\xi_n}p(\xi',\xi_n)}
{|\xi'|^2+\xi^2_n}\right]&=&\frac{\partial_{\xi_n}p(\xi',\xi_n)|
_{\xi_n=i|\xi'|}}{2i|\xi'|(\xi_n-i|\xi'|)};\\
\pi^+_{\xi_n}\left[\frac{2\xi_np(\xi',\xi_n)}{(|\xi'|^2+\xi^2_n)^2}\right]
&=&\frac{p(\xi',i|\xi'|)}
{2i|\xi'|(\xi_n-i|\xi'|)^2}+\frac{\partial_{\xi_n}p(\xi',\xi_n)|
_{\xi_n=i|\xi'|}}{2i|\xi'|(\xi_n-i|\xi'|)},
\end{eqnarray*}
\noindent so,
$$\pi^+_{\xi_n}\partial_{\xi_n}\sigma_L(F)(\xi',\xi_n)=\frac{-p(\xi',i|\xi'|)}
{2i|\xi'|(\xi_n-i|\xi'|)^2}.\eqno(7.2)$$ \noindent Using (7.2),
then
\\
$${\rm trace}[\pi^+_{\xi_n}\partial_{\xi_n}\sigma_L(F)(\xi',\xi_n)\times
 \partial_{\xi_n}^{2}\sigma_L(F)(\xi',\xi_n)]$$
$$
={\rm trace}\left\{\frac{-p(\xi',i|\xi'|)}
{2i|\xi'|(\xi_n-i|\xi'|)^2}\times\left[\frac{\partial^2_{\xi_n}p(\xi',\xi_n)}{|\xi'|^2+\xi_n^2}-
\frac{4\xi_n\partial_{\xi_n}p(\xi',\xi_n)}{(|\xi'|^2+\xi_n^2)^2}-\frac{2|\xi'|^2-6\xi_n^2}{(|\xi'|^2+\xi_n^2)^3}
p(\xi',\xi_n)\right ]\right\}$$
$$=ia_{n,m}\left[\frac{1}{(1+i\xi_n)^2(1+\xi_n^2)}+\frac{4i\xi_n}{(1+\xi_n^2)^2(1+i\xi_n)}+\frac{1-3\xi_n^2}
{(1+\xi_n)^3}\right].$$

\noindent Integrate with respect to $\xi_n$, then\\
$$ia_{n,m}\int_{\Gamma+}
\left[\frac{1}{(1+i\xi_n)^2(1+\xi_n^2)}+\frac{4i\xi_n}
{(1+\xi_n^2)^2(1+i\xi_n)}+\frac{1-3\xi_n^2}{(1+\xi_n)^3}\right ]d\xi_n\\
=3\pi i.$$
\noindent So,
\begin{eqnarray*}
D_4&=&(-i)^3\partial_{x_n} f_1(x',0)\partial_{x_n}f_2(x',0)
\int_{|\xi'|=1}
 3\pi i\sigma(\xi')dx_1\wedge
 dx_2\\
&=&-6{\pi}^2\partial_{x_n}
f_1(x',0)\partial_{x_n}f_2(x',0)dx_1\wedge
 dx_2\\
&=&-6{\pi}^2\partial_{x_n} f_1(x',0)\partial_{x_n}f_2(x',0){\rm
Vol}_{Y}.
\end{eqnarray*}
\noindent Then we deduce the formula:
$$\Omega_{2,{\rm flat}}(f_1,f_2)=D_1+D_4=-3{\pi}^2d(f_1|_{Y})\wedge\star d(f_2|_{Y})
-6{\pi}^2\partial_{x_n} f_1(x',0)\partial_{x_n}f_2(x',0){\rm
Vol}_{Y}.\eqno(7.3)$$

\noindent{\bf{Remark:}}~~$X$ has the product structure near the
boundary, so $\partial_{x_n}f_i|_{x_n=0}\in
C^{\infty}(Y)$.\\
\noindent{\bf c)~~Non-flat Case}\\
\indent Since $n=3$ and $r\leq-1,s\leq-1$,  so $~r=s=-1,~ |\alpha|=k=j=0$
and $|\beta|=|\delta|=1$, by (3.19) we have:\\
$$\Omega_{2}(f_1,f_2)=\sum_{|\beta|=1}\sum_{|\delta|=1}(-i)^3\partial_{x}^{\beta}
f_1(x',0)\partial_{x}^{\delta}f_2(x',0)\times$$
 $$ \int_{|\xi'|=1}\int^
 {+\infty}_{-\infty}
 {\rm trace}\left[\pi^+_{\xi_n}\partial^{\beta}_{\xi}\sigma_L(F)(\xi',\xi_n)\times
 \partial_{\xi_n}\partial^{\delta}_{\xi}\sigma_L(F)(\xi',\xi_n)\right]d\xi_n\sigma(\xi')dx_1\wedge
 dx_2.\eqno(7.4)$$
 \indent Observe: (7.4) is similar to case b) and the only
 difference is that $\sigma_L(F)$ is dependent of $x$. In the spirit
 of [10], we compute this form by the normal coordinate
 way.\\
 \indent In (7.4), there is no
 derivative $\partial_{x_i}$ with respect to {\rm trace}, so we take the
 normal coordinate and take boundary point $x=x_0$. Then
 $g^{ij}(x_0)=\delta_{ij}$ where $[g^{ij}]$ is the inverse matrix
 of metric matrix and this case is same to the case b). Whereas:
\begin{eqnarray*}
df_1(x',0)\wedge\star df_2(x',0)|_{x_0}
&=&\sum^{2}_{i,j=1}\partial_{x_i}f_1(x_0)\partial_{x_j}f_2(x_0)g^{ij}(x_0){\rm
det}^{\frac{1}{2}}{[g_{ij}(x_0)]}dx_1
\wedge dx_2\\
&=&\sum^{2}_{i=1}\partial_{x_i}f_1(x_0)\partial_{x_i}f_2(x_0)dx_1
\wedge dx_2
\end{eqnarray*}
and ${\rm Vol}_{Y}|_{x_0}=dx_1\wedge dx_2,$
so (7.3) is correct in this case. By [3],\\
$$\Omega_{2}(f_1|_Y,f_2|_Y)=-8\pi d(f_1|_Y)\wedge\star
d(f_2|_Y),$$ \noindent then we obtain:\\
\noindent{\bf Theorem 7.1}~~Suppose that $(X,g)$ is a
3-dimensional compact oriented Riemannian manifold with boundary
$Y$ and has product metric near the boundary, then we have: {\it
\begin{eqnarray*} \Omega_{2}(f_1,f_2) &=&-3{\pi}^2
d(f_1|_Y)\wedge\star d(f_2|Y)-6{\pi}^2\partial_{x_n}
f_1(x',0)\partial_{x_n}f_2(x',0){\rm Vol}_{Y}\\
&=&\frac{3}{8}\pi\Omega_{2}(f_1|_Y,f_2|_Y)-6{\pi}^2\partial_{x_n}
f_1(x',0)\partial_{x_n}f_2(x',0){\rm Vol}_{Y},
\end{eqnarray*}
$~~~~~~\widetilde{{\rm
Wres}}(\pi^+f_0[\pi^{+}F,\pi^+f_1][\pi^+F,\pi^+f_2])$
$$= -3{\pi}^2\int_Yf_0|_Y\left[d(f_1|_Y)\wedge\star
d(f_2|Y)+2\partial_{x_n} f_1(x',0)\partial_{x_n}f_2(x',0){\rm
Vol}_{Y}\right].\eqno(7.5)$$ } \indent The above formula is the
generalization to  manifolds with boundary of idea in [3] when
$n=3$. By [3], $\Omega_{2}(f_1|_Y,f_2|_Y)$ is conformally
invariant. So although $\Omega_{2}(f_1,f_2)$ is not a conformal
invariant,
but we have:\\
\noindent{\bf Corollary 7.2}~~{\it
$$\Omega_{2}(f_1,f_2)+6{\pi}^2\partial_{x_n}
f_1(x',0)\partial_{x_n}f_2(x',0){\rm Vol}_Y$$ is a conformal
invariant of $(X,g)$.\\}
 \noindent{\bf Acknowledgement:}~~The
author would like to thank Professors Huitao Feng and Weiping
Zhang for their helpful discussions. He also thanks the referee
for his
careful reading and helpful comments.\\

{\footnotesize

\end{document}